# Path following control for a reversing general 2-trailer system

Oskar Ljungqvist[1], Daniel Axehill[1], Anders Helmersson[1]

*Abstract*— In order to meet the requirements for autonomous systems in real world applications, reliable path following controllers have to be designed to execute planned paths despite the existence of disturbances and model errors. In this paper we propose a Linear Quadratic controller for stabilizing a 2-trailer system with possible off-axle hitching around preplanned paths in backward motion. The controller design is based on a kinematic model of a general 2-trailer system including the possibility for off-axle hitching. Closed-loop stability is proved around a set of paths, typically chosen to represent the possible output from the path planner, using theory from linear differential inclusions. Using convex optimization tools a single quadratic Lyapunov function is computed for the entire set of paths.

## I. INTRODUCTION

The past decades have witnessed a rapid increase in autonomous vehicle research. Historically Advanced Driver Assistance Systems (ADSA) have been designed to enhance safety like Lane Keeping Assist (LKA), Adaptive Cruise control (ACC) and Automatic Breaking. In recent years ADSA systems have been developed to help the driver in more complex tasks, such as reversing with a trailer [1]. A fair amount of training is necessary for an unexperienced driver to smoothly perform such a task. Adding one more degree of freedom to the problem, by considering a truck with a dolly steered trailer, makes it almost impossible to reverse for a nonprofessional driver. An unexperienced driver will encounter problems already when performing simple tasks, such as reversing in a straight line. In this paper we propose a Linear Quadratic (LQ) controller with feed forward action to stabilize the lateral and angular error dynamics for a truck with a dolly steered trailer around a precomputed path in backward motion. The closed-loop system is analyzed and closed-loop stability is proven for a set of paths to guarantee safe execution.

### A. Related work

The nonlinear dynamics of the standard n-trailer, which assume on-axle hitching connections, are well explored and derivation of the nonlinear kinematic model can be found in [2]. The standard n-trailer possesses the property that it is differentially flat when using the position of the rear axle of the last trailer as the flat outputs [3]. This flatness property is used in [4] to design controllers based on feedback linearization. In [2] and [5] it is shown that it is possible, after a change of coordinates, to put the system into a chained form which makes it very easy to perform motion planning in an obstacle free environment. However, the on-axle hitching assumption does not hold in general for regular truck and trailer systems nor for the truck that will be used in this work. The nonlinear dynamics of the general n-trailer system, which allow off-axle hitching, are derived in [6]. As shown in [3], the general 1-trailer is still differentially flat using a certain choice of flat outputs. However, the flatness property does not hold for the general n-trailer system ($n \geq 2$). In [7] this problem is circumvented by introducing a "ghost vehicle" which has equivalent stationary behavior but different transient behavior. In [8] they instead design a controller using input-output linearization to stabilize the general n-trailer system around a path of constant curvature. This approach is limited to forward motion since the introduced zero dynamics is unstable in backward motion. In [9], they design an LQ-controller to stabilize the general 2-trailer system around different equilibrium configurations corresponding to straight lines and circles in forward and backward motion. A cascade control approach to stabilize the general 2-trailer system around piecewise linear reference paths is proposed in [10]. This approach can be motivated off-line in a motion planner but when on-line execution of a precalculated path in considered this method is not using all information about the path. A complete survey in the area of control of truck and trailer systems can be found in [11]. A similar approach as in [9] is used in this paper, but instead of restricting to stabilize the general 2-trailer system around circles of different radiuses we propose a framework that can be proven to stabilize the system around any precalculated path from a specified set of possible paths given by a motion planner [12]. The advantage of this approach is that the path is obtained from the true system dynamics making it, in theory, kinematically feasible to follow perfectly. To facilitate the stability analysis of nonlinear systems one possible technique is to describe them as Linear Differential Inclusions (LDIs) [13]. In [14] they use Linear Matrix Inequalities (LMIs) techniques together with LDI techniques to show stability for a gas turbine engine by finding a common Lyapunov function by approximating the convex hull of the derivative of the closed-loop system matrix. An inner approximation of the convex hull is obtained by differentiating the closed-loop system around some typical working points. A similar approach will be used in this work but instead of doing an inner approximation of the convex hull we do an outer approximation using global optimization tools [15].

The paper is organized as follows: in Section II the lateral and angular error dynamics of the general 2-trailer system is modeled using the Frenet frame; in Section III an LQ-controller is proposed to stabilize the lateral and angular error dynamics in backward motion; in Section IV a sufficient condition to guarantee local asymptotic stability around a set of paths is proven and in Section V the method is applied to a particular problem in simulations.

[1]Division of Automatic Control, Linköping University, Sweden, (e-mail: {oskar.ljungqvist, daniel.axehill, anders.helmersson}@liu.se)

## II. MODELING

In this section we derive a model for the error dynamics around a path for a general 2-trailer system. A schematic description of the system is presented in Fig. 1. The system has an off-axle connection between the truck and the dolly and an on-axle connection between the dolly and the trailer. With the generalized coordinates $\mathbf{p} = [x_3, y_3, \theta_3, \beta_3, \beta_2]^T$ where $(x_3, y_3)^T$ denotes the global position of the center of the rear axle of the trailer, $\theta_3$ denotes the global orientation of the trailer, $\beta_3$ is the relative angle between the trailer and the dolly and $\beta_2$ is the relative angle between the dolly and the truck. The geometric length $L_3$ denotes the distance between the rear axle of the trailer and the axle of the dolly, $L_2$ denotes the distance between the axle of the dolly and the off-axle hitching connection at the truck, $M_1$ is the length of the off-axle hitching and $L_1$ is the distance between the rear axle of the truck and its front axle. The steering angle $\alpha$ is a control signal and $v_3$ is defined as the velocity of the trailer. The sign of $v_3$ decides the direction of motion where $v_3 > 0$ denotes forward motion and $v_3 < 0$ denotes backward motion. A recursive formula for the general n-trailer based on nonholonomic and holonomic constraints is derived in [6] and a model for the specific 2-trailer case with off-axle hitching at the truck is presented in [9]:

$$\dot{x}_3 = v_3 \cos\theta_3 \tag{1a}$$
$$\dot{y}_3 = v_3 \sin\theta_3 \tag{1b}$$
$$\dot{\theta}_3 = v_3 \frac{\tan\beta_3}{L_3} \tag{1c}$$
$$\dot{\beta}_3 = v_3 \left( \frac{\tan\beta_2 - \frac{M_1}{L_1}\tan\alpha}{L_2 \cos\beta_3 C_1(\beta_2, \alpha)} - \frac{\tan\beta_3}{L_3} \right) \tag{1d}$$
$$\dot{\beta}_2 = v_3 \left( \frac{\frac{\tan\alpha}{L_1 \cos\beta_2} - \frac{\tan\beta_2}{L_2} + \frac{M_1}{L_1 L_2}\tan\alpha}{\cos\beta_3 C_1(\beta_2, \alpha)} \right) \tag{1e}$$

where

$$C_1(\beta_2, \alpha) = \left(1 + \frac{M_1}{L_1}\tan\beta_2 \tan\alpha\right) \tag{2}$$

Represent (1) as $\dot{\mathbf{p}} = v_3 f(\mathbf{p}, \alpha)$. This model is valid under low speed assumptions since it is derived based on holonomic and nonholonomic constraints where a no slip assumption has to hold. The system is assumed to operate on a relatively flat surface. A method known as time scaling [16] can be applied to eliminate the speed dependence in the model and therefore we assume without loss of generality $v_3$ to be an a priori decided constant. The model (1) has singularities which is discussed in [6]. In order to avoid them we must bound the relative angles as $|\beta_3| < \pi/2$ and $|\beta_2| < \pi/2$. From (1a)-(1e) we can conclude that the steering angle $\alpha$ enters as $\tan\alpha$ everywhere and therefore we do the input substitution

$$u = \tan\alpha \tag{3}$$

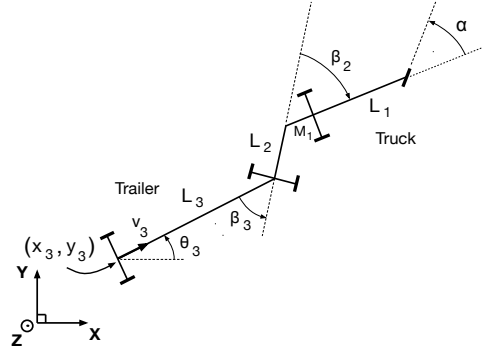

Fig. 1: Definition of the geometric lengths and the generalized coordinates used to model the general 2-trailer system.

and for what follows we assume that we are manipulating $u$ instead of $\alpha$. Since $\tan\alpha \in ]-\infty, \infty[$ for $\alpha \in ]-\pi/2, \pi/2[$ we do not need to constrain $u$ which is allowed to take all real values. In practice, we have saturation in the steering angle which has to be taken into account when designing controller such that it does not saturate the steering angle.

### A. Frenet frame

Given a feasible trajectory $(\mathbf{p}_0(\cdot), u_0(\cdot))$ satisfying

$$\dot{\mathbf{p}}_0 = v_3 f(\mathbf{p}_0(t), u_0(t)) \tag{4}$$

we can describe the system (1a)-(1e) in terms of the deviation from a trajectory generated by the system in (4), see Fig. 2. Introduce the parameter $s$ as the distance traveled by the rear axle of the trailer onto its projection to the path up to time $t$. With this new parameterization $s$ equation (4) can be written as

$$\frac{d\mathbf{p}_0}{ds} = f(\mathbf{p}_0(s), u_0(s)) \tag{5}$$

Using standard geometry the curvature $\kappa_0(s)$ of the path followed by the trailer is given by

$$\kappa_0(s) = \frac{d\theta_{3,0}}{ds} = \frac{\tan\beta_{3,0}(s)}{L_3} \tag{6}$$

Denote $z_3$ as the lateral deviation between the midpoint of the rear axle of the trailer onto its projection to the path. Let $\tilde{\theta}_3(t) = \theta_3(t) - \theta_{3,0}(s(t))$ be the error in the heading of the trailer with respect to the path, define $\tilde{\beta}_3(t) = \beta_3(t) - \beta_{3,0}(s(t))$, $\tilde{\beta}_2(t) = \beta_2(t) - \beta_{2,0}(s(t))$ and $\tilde{u}(t) = u(t) - u_0(s(t))$. Then system (1a)-(1e) can be described in the

Frenet frame, using the chain rule, as

$$\dot{s} = v_3 \frac{\cos \tilde{\theta}_3}{1 - \kappa_0(s)z_3} \quad (7a)$$

$$\dot{z}_3 = v_3 \sin \tilde{\theta}_3 \quad (7b)$$

$$\dot{\tilde{\theta}}_3 = v_3 \left( \frac{\tan(\tilde{\beta}_3 + \beta_{3,0})}{L_3} - \frac{\kappa_0(s)\cos \tilde{\theta}_3}{1 - \kappa_0(s)z_3} \right) \quad (7c)$$

$$\dot{\tilde{\beta}}_3 = v_3 \left( \frac{\tan(\tilde{\beta}_2 + \beta_{2,0}) - \frac{M_1}{L_1}(\tilde{u} + u_0)}{L_2 \cos(\tilde{\beta}_3 + \beta_{3,0}) C_2(\tilde{\beta}_2 + \beta_{2,0}, \tilde{u} + u_0)} - \frac{\tan(\tilde{\beta}_3 + \beta_{3,0})}{L_3} \right.$$
$$\left. - \frac{\cos \tilde{\theta}_3}{1 - \kappa_0(s)z_3} \left( \frac{\tan \beta_{2,0} - \frac{M_1}{L_1}u_0}{L_2 \cos \beta_{3,0} C_2(\beta_{2,0}, u_0)} - \kappa_0(s) \right) \right) \quad (7d)$$

$$\dot{\tilde{\beta}}_2 = v_3 \left( \left( \frac{\frac{\tilde{u}+u_0}{L_1 \cos(\tilde{\beta}_2 + \beta_{2,0})} - \frac{\tan(\tilde{\beta}_2+\beta_{2,0})}{L_2} + \frac{M_1}{L_1 L_2}(\tilde{u}+u_0)}{\cos \beta_3 C_2(\tilde{\beta}_2 + \beta_{2,0}, \tilde{u}+u_0)} \right) \right.$$
$$\left. - \frac{\cos \tilde{\theta}_3}{1 - \kappa_0(s)z_3} \left( \frac{\frac{u_0}{L_1 \cos \beta_{2,0}} - \frac{\tan \beta_{2,0}}{L_2} + \frac{M_1}{L_1 L_2} u_0}{\cos \beta_{3,0} C_2(\beta_{2,0}, u_0)} \right) \right) \quad (7e)$$

where

$$C_2(\beta_2, u) = 1 + \frac{M_1}{L_1} \tan \beta_2 u \quad (8)$$

This transformation is valid in a tube around the desired path for which the denominator $1 - \kappa_0(s)z_3$ is greater than zero. The width of this tube depends on the curvature $\kappa_0(s)$ and when the curvature tends to zero (a straight line) $z_3$ can vary arbitrarily. Essentially to avoid singularities in the transformation we must have that $|z_3| < |\kappa_0(s)^{-1}|$ when $z_3$ and $\kappa_0(s)$ have the same sign. To guarantee $\dot{s}$ to be a monotonic function in time we assume $|\tilde{\theta}_3| < \pi/2$. In this paper we are not interested in how fast the nominal trajectory (4) is covered but only that the resulting path is traversed with low lateral and angular tracking error. Since we are only interested in the lateral and the angular dynamics, the equation in (7a) is disregarded and the dependence on $s$ in (7b)-(7e) will be treated as a parametric uncertainty. Since the longitudinal equation (7a) has been disregarded the nominal trajectory generated by (4) will be considered as a nominal path. From (7b)-(7e) we can conclude that the tracking error dynamics do not depend explicitly on $x_{3,0}$, $y_{3,0}$ nor $\theta_{3,0}$ and therefore we introduce a new state vector $\mathbf{p}_e = [z_3, \tilde{\theta}_3, \beta_3, \beta_2]^T$. With $\mathbf{p}_e$ as state vector the coordinates for the nominal path is $\mathbf{p}_{e,0} = [0, 0, \beta_{3,0}, \beta_{2,0}]^T$ and the corresponding error state vector is $\tilde{\mathbf{p}}_e = \mathbf{p}_e - \mathbf{p}_{e,0} = [z_3, \tilde{\theta}_3, \tilde{\beta}_3, \tilde{\beta}_2]^T$. The nonlinear tracking error dynamics can be written as

$$\dot{\tilde{\mathbf{p}}}_e = f_e(\tilde{\mathbf{p}}_e, \tilde{u}, \mathbf{p}_{e,0}, u_0) \quad (9)$$

Furthermore, we focus only on stabilizing the nonlinear system (9) in some neighborhood around the nominal path, i.e., $(\tilde{\mathbf{p}}_e, \tilde{u}) = (\bar{0}, 0)$. By construction the point $(\tilde{\mathbf{p}}_e, \tilde{u}) = (\bar{0}, 0)$ is an equilibrium to (9) since $f_e(\bar{0}, 0, \mathbf{p}_{e,0}, u_0) = \bar{0}$.

### B. Linearization around paths

The nonlinear system (7b)-(7e) can be linearized around the nominal path $(\mathbf{p}_{e,0}, u_0)$ by linearizing (9) around $(\tilde{\mathbf{p}}_e, \tilde{u}) =$

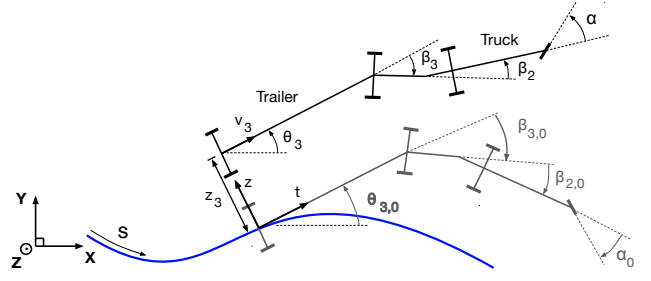

Fig. 2: An illustrative description of the Frenet frame with its moving coordinate system located at the projection of the rear axle of the trailer to the blue path $(x_{3,0}, y_{3,0})(s)$. The black truck and trailer is the actual system and the gray truck and trailer is the reference vehicle or the desired vehicle configuration at this specific value of $s(t)$.

$(\bar{0}, 0)$. This point is an equilibrium to (9) and hence a first order Taylor series expansion yields

$$\dot{\tilde{\mathbf{p}}}_e = A(\bar{0}, 0, \mathbf{p}_{e,0}, u_0)\tilde{\mathbf{p}}_e + B(\bar{0}, 0, \mathbf{p}_{e,0}, u_0)\tilde{u} \quad (10)$$

where

$$A(\tilde{\mathbf{p}}_e, \tilde{u}, \mathbf{p}_{e,0}, u_0) = \frac{\partial f_e}{\partial \tilde{\mathbf{p}}_e}(\tilde{\mathbf{p}}_e, \tilde{u}, \mathbf{p}_{e,0}, u_0) \quad (11)$$

$$B(\tilde{\mathbf{p}}_e, \tilde{u}, u_0, \mathbf{p}_{e,0}) = \frac{\partial f_e}{\partial \tilde{u}}(\tilde{\mathbf{p}}_e, \tilde{u}, \mathbf{p}_{e,0}, u_0) \quad (12)$$

The matrices $A(\tilde{\mathbf{p}}_e, \tilde{u}, \mathbf{p}_{e,0}, u_0)$ and $B(\tilde{\mathbf{p}}_e, \tilde{u}, \mathbf{p}_{e,0}, u_0)$ can be found in the Appendix. For control design we intend to use an LQ-controller to stabilize (9) around the origin $(\tilde{\mathbf{p}}_e, \tilde{u}) = (\bar{0}, 0)$ and for control design we are thus only interested in the linearization around the nominal path $(\mathbf{p}_{e,0}, u_0) = (\bar{0}, 0)$ (a straight line) and the relevant matrices $A$ and $B$ are then

$$A = v_3 \begin{bmatrix} 0 & 1 & 0 & 0 \\ 0 & 0 & \frac{1}{L_3} & 0 \\ 0 & 0 & -\frac{1}{L_3} & \frac{1}{L_2} \\ 0 & 0 & 0 & -\frac{1}{L_2} \end{bmatrix} \quad B = v_3 \begin{bmatrix} 0 \\ 0 \\ -\frac{M_1}{L_1 L_2} \\ \frac{L_2 + M_1}{L_1 L_2} \end{bmatrix} \quad (13)$$

The characteristic polynomial for the matrix $A$ is thus given by

$$\det(\lambda I - A) = \lambda^2 \left(\lambda + \frac{v_3}{L_3}\right)\left(\lambda + \frac{v_3}{L_2}\right) \quad (14)$$

from which it follows that the nonlinear system is, as one would expect, unstable in backward motion ($v_3 < 0$) since the linear system has two poles in the right half plane corresponding to the two hitching connections of the vehicle configurations. In forward motion ($v_3 > 0$) the system is marginally stable because of the double integrator. Due to the off-axle hitching the system has a zero in some of the output channels, as concluded in [9]. This zero is stable in backward motion and unstable in forward motion.

### III. STABILIZATION

For the basics in Linear Quadratic (LQ) control we refer to [17]. To stabilize the nonlinear system (9) we design an

LQ-controller
$$\tilde{u} = -K\tilde{\mathbf{p}}_e \quad (15)$$

based on the linearized model (10) around the nominal trajectory $(\mathbf{p}_{e,0}, u_0) = (\bar{0}, 0)$, i.e., a straight line. As long as the velocity $v_3 \neq 0$ the nonlinear system (9) is locally controllable around the nominal path $(\mathbf{p}_{e,0}, u_0) = (\bar{0}, 0)$ with linear state feedback since the pair $(A, B)$ is controllable. According to standard LQ theory, the feedback gain $K$ is chosen such that the cost function

$$\mathscr{J} = \int_0^\infty \left( \tilde{x}^T Q \tilde{x} + \tilde{u}^2 \right) dt \quad (16)$$

is minimized w.r.t $\tilde{u}$ where $Q = Q^T \succeq 0$. The optimal feedback gain is

$$K = -B^T P \quad (17)$$

where $P = P^T \succ 0$ is the unique solution to the Algebraic Riccati Equation

$$A^T P + PA + Q = PBB^T P \quad (18)$$

It can be shown that the feedback gain $K$ is actually independent of the size of $v_3$ and the locations of the poles to the closed-loop system are independent of the sign of $v_3$. This is of big importance and we summarize it in the following theorem.

*Theorem 1:* Consider LQ-control of the linear system

$$\dot{x} = v(Ax + Bu) \quad (19)$$

with the weight matrices $Q = C^T C \succeq 0$ and $R = R^T \succ 0$. Assume $v \neq 0$ is a constant, $(A, B)$ is controllable and $(C, A)$ is observable. Then
1) the optimal feedback gain $K$ is independent of the magnitude of v.
2) the closed-loop system $\dot{x} = v(A - BK)x$ has the same pole locations for $v$ and $-v$ when using LQ-control.

*Proof:* 1): Assume that $v = a\tilde{v}$ where $a$ is a positive constant. The ARE for (19) becomes

$$A^T \tilde{v}aP + \tilde{v}aPA + Q = aP\tilde{v}BR^{-1}B^T \tilde{v}aP$$

This equation can be rewritten as

$$A^T \tilde{v}\tilde{P} + \tilde{v}\tilde{P}A + Q = \tilde{P}\tilde{v}BR^{-1}B^T \tilde{v}\tilde{P}$$

where we have defined $\tilde{P} = Pa$. The feedback gain becomes

$$K = -R^{-1}\tilde{v}B^T aP = -R^{-1}\tilde{v}B^T \tilde{P}$$

which is independent of $a$.
2): Introduce the Hamiltonian matrix:

$$H(v) = \begin{bmatrix} vA & -v^2 BR^{-1}B^T \\ -Q & -vA^T \end{bmatrix}$$

This matrix has the property that it has all its eigenvalues symmetric w.r.t the real and imaginary axis. Using LQ-control the closed-loop system poles are the negative eigenvalues to $H$. Define

$$T = \begin{bmatrix} \mathbb{I} & 0 \\ 0 & -\mathbb{I} \end{bmatrix} \quad (20)$$

and apply the similarity transformation $H \mapsto THT^{-1}$. The result is

$$THT^{-1} = \begin{bmatrix} vA & v^2 BR^{-1}B^T \\ Q & -vA^T \end{bmatrix} = -H(-v) \quad (21)$$

Since the eigenvalues to $H$ and $THT^{-1}$ coincide and because the eigenvalues to $H$ and $-H$ also coincide we can conclude that the eigenvalues to H must be the same for $v$ and $-v$. In particular this holds for the negative eigenvalues, which are the location of the poles for the closed-loop system. ■

*Remark 1:* Even though the poles have the same locations we can not tell anything about where the zeros are located.

In practice Theorem 1 says that it is, without loss of generality, valid to assume that the velocity of the trailer, $v_3$, only takes on the values in the index set $\mathscr{I} = \{-1, 1\}$ when designing linear state feedback controllers for the system (10). Alternatively, similar to [9] a method known as time-scaling [16] can be applied to eliminate the velocity dependence.

## IV. STABILITY ANALYSIS

When applying the feedback controller (15) to the nonlinear system (9) the autonomous closed-loop system can be written in a compact form as

$$\dot{\tilde{\mathbf{p}}}_e = f_{cl}(\tilde{\mathbf{p}}_e, \mathbf{p}_{e,0}, u_0) \quad (22)$$

where $\tilde{\mathbf{p}}_e = \bar{0}$ is an equilibrium, since $f_{cl}(\bar{0}, \mathbf{p}_{e,0}, u_0) = \bar{0}$. It will not be possible to show global asymptotic stability of (22) since the Frenet frame description is not globally well-defined, as discussed in Section II. However it is possible to show local asymptotic stability around a set of nominal paths satisfying (4). The idea in this work is to, e.g., choose this set to represent the possible output from a path planner as the one in [12].

*Assumption 1:* The set of paths is given by $(\mathbf{p}_{e,0}, u_0) \in \mathbb{P} \times \mathbb{U}$ where $\mathbb{P}$ and $\mathbb{U}$ are assumed bounded. Assume $f_{cl} : \tilde{\mathbb{P}} \times \mathbb{P} \times \mathbb{U} \to \mathbb{R}^4$ is continuously differentiable w.r.t. $\tilde{\mathbf{p}}_e \in \tilde{\mathbb{P}} = \{\tilde{\mathbf{p}}_e \in \mathbb{R}^4 \mid \|\tilde{\mathbf{p}}_e\|_2 < r\}$ and the Jacobian matrix $[\partial f_{cl}/\partial \tilde{\mathbf{p}}_e]$ is bounded and Lipschitz on $\tilde{\mathbb{P}}$, uniformly in $(\mathbf{p}_{e,0}, u_0)$.

*Theorem 2:* Consider the closed-loop system (22). Under Assumption 1, let

$$A_{cl}(\mathbf{p}_{e,0}, u_0) = \frac{\partial f_{cl}}{\partial \tilde{\mathbf{p}}_e}(\bar{0}, \mathbf{p}_{e,0}, u_0) \quad (23)$$

and define the convex hull

$$\bar{\mathbb{A}}_{cl} = \mathbf{Co} \; \{A_{cl}(\mathbf{p}_{e,0}, u_0) \in \mathbb{R}^{4 \times 4} \mid (\mathbf{p}_{e,0}, u_0) \in \mathbb{P} \times \mathbb{U}\} \quad (24)$$

If there exist matrices $P = P^T \succ 0$ and $Q = Q^T \succ 0$ that satisfy

$$PA_{cl} + A_{cl}^T P \preceq -Q \quad \forall A_{cl} \in \bar{\mathbb{A}}_{cl} \quad (25)$$

then $\tilde{\mathbf{p}}_e = \bar{0}$ is an exponentially stable equilibrium point to (22) for all $(\mathbf{p}_{e,0}, u_0) \in \mathbb{P} \times \mathbb{U}$.

*Proof:* The proof is based on results from [13] and [18]. Introduce the quadratic Lyapunov function candidate $V(\tilde{\mathbf{p}}_e) = \tilde{\mathbf{p}}_e^T P \tilde{\mathbf{p}}_e$. For an arbitrary pair $(\mathbf{p}_{e,0}, u_0) \in \mathbb{P} \times \mathbb{U}$, Assumption 1 guarantees existence of a function $g(\tilde{\mathbf{p}}_e, \mathbf{p}_{e,0}, u_0)$ that, for some $L > 0$, satisfies

$$\|g(\tilde{\mathbf{p}}_e, \mathbf{p}_{e,0}, u_0)\|_2 < L\|\tilde{\mathbf{p}}_e\|_2^2$$

and of positive constants $c_1$ and $c_2$ such that

$$\dot{V}(\tilde{\mathbf{p}}_e) = \tilde{\mathbf{p}}_e^T P f_{cl}(\tilde{\mathbf{p}}_e, \mathbf{p}_{e,0}, u_0) + f_{cl}^T(\tilde{\mathbf{p}}_e, \mathbf{p}_{e,0}, u_0) P \tilde{\mathbf{p}}_e$$
$$= \tilde{\mathbf{p}}_e^T (PA_{cl}(\mathbf{p}_{e,0}, u_0) + PA_{cl}^T(\mathbf{p}_{e,0}, u_0))) \tilde{\mathbf{p}}_e + 2\tilde{\mathbf{p}}_e^T P g(\tilde{\mathbf{p}}_e, \mathbf{p}_{e,0}, u_0)$$
$$\leq -\tilde{\mathbf{p}}_e^T Q \tilde{\mathbf{p}}_e + 2\tilde{\mathbf{p}}_e^T P g(\tilde{\mathbf{p}}_e, \mathbf{p}_{e,0}, u_0)$$
$$\leq -(c_1 - 2c_2 L \rho) \|\tilde{\mathbf{p}}_e\|_2^2, \quad \forall \|\tilde{\mathbf{p}}_e\|_2 < \rho$$

If we choose $\rho = \min\{r, \frac{c_1}{2c_2 L}\}$ it holds that $\dot{V} < 0$ in $\|\tilde{\mathbf{p}}_e\|_2 < \rho$. Furthermore if we can find a common $P = P^T \succ 0$ and $Q = Q^T \succ 0$ that satisfies

$$PA_{cl}(\mathbf{p}_{e,0}, u_0) + PA_{cl}^T(\mathbf{p}_{e,0}, u_0) \preceq -Q \quad \forall (\mathbf{p}_{e,0}, u_0) \in \mathbb{P} \times \mathbb{U} \quad (26)$$

then the origin is an exponentially stable equilibrium for all parameters satisfying $(\mathbf{p}_{e,0}, u_0) \in \mathbb{P} \times \mathbb{U}$. Finally, using results from [13], (25) implies (26) and the proof is complete. ∎

The convex hull $\bar{\mathbb{A}}_{cl}$ is generally hard to find [13] and in this paper we have chosen to do an outer approximation of the convex hull $\bar{\mathbb{A}}_{cl}$. The differentiation of the open-loop system (9) can be found in the Appendix. Denoting the feedback $K = \begin{bmatrix} k_1 & k_2 & k_3 & k_4 \end{bmatrix}$ and using the index notation in the Appendix for the index values in $A(\bar{0}, 0, \mathbf{p}_{e,0}, u_0)$ and $B(\bar{0}, 0, \mathbf{p}_{e,0}, u_0)$ the matrix $A_{cl}$ can be expressed as

$$A_{cl} = v_3 \begin{bmatrix} 0 & 1 & 0 & 0 \\ a_{21} & 0 & a_{23} & 0 \\ a_{31} - b_3 k_1 & -b_3 k_2 & a_{33} - b_3 k_3 & a_{34} - b_3 k_4 \\ a_{41} - b_4 k_1 & -b_4 k_2 & a_{43} - b_4 k_3 & a_{44} - b_4 k_4 \end{bmatrix} \quad (27)$$

Each value in this matrix can be bounded as

$$A_{cl}(i,j)_{\min} \leq A_{cl}(i,j) \leq A_{cl}(i,j)_{\max} \quad (28)$$

where $A_{cl}(i,j)_{\min}$ and $A_{cl}(i,j)_{\max}$ are defined as

$$A_{cl}(i,j)_{\min} = \min_{(\mathbf{p}_{e,0}, u_0) \in \mathbb{P} \times \mathbb{U}} A_{cl}(i,j)(\mathbf{p}_{e,0}, u_0) \quad (29a)$$

$$A_{cl}(i,j)_{\max} = \max_{(\mathbf{p}_{e,0}, u_0) \in \mathbb{P} \times \mathbb{U}} A_{cl}(i,j)(\mathbf{p}_{e,0}, u_0) \quad (29b)$$

Based on (29a) and (29b), we can outer bound $\bar{\mathbb{A}}_{cl}$ in a 10 dimensional box. From the theory for polytopic LDIs in [13], it follows that it suffices to consider the $2^{10} = 1024$ vertices of the box. Hence a sufficient condition for (25) to hold is that there exists a common matrix $P = P^T \succ 0$ and $Q = Q^T \succ 0$ such that

$$A_{cl,i}^T P + PA_{cl,i} \preceq -Q \quad i = 1, 2, ..., 1024 \quad (30)$$

If such a $P$ exists then Theorem 2 guarantees local exponential stability for the closed-loop system (22) around the origin for any $(\mathbf{p}_{e,0}, u_0) \in \mathbb{P} \times \mathbb{U}$. To ensure a certain decay rate, i.e. to ensure $dV(\tilde{\mathbf{p}}_e)/dt \leq -2\varepsilon V(\tilde{\mathbf{p}}_e)$, $\varepsilon > 0$ of the polytopic LDI

$$\dot{\tilde{\mathbf{p}}}_e = A_{cl} \tilde{\mathbf{p}}_e, \quad A_{cl} \in \mathbf{Co}\{A_{cl,1}, ..., A_{cl,1024}\} \quad (31)$$

we can replace condition (30) with

$$A_{cl,i}^T P + PA_{cl,i} + 2\varepsilon P \preceq 0 \quad i = 1, 2, ..., 1024 \quad (32)$$

This will ensure a decay rate of at least $\varepsilon$ for (22) [13] around the origin.

## V. RESULTS

In this section we start by applying Theorem 2 and the strategy described in Section IV on a general 2-trailer system in backward motion ($v_3 < 0$) with the geometric lengths $L_1 = 3.8$ m, $L_2 = 2.8$ m, $M_1 = 0.72$ m and $L_3 = 6.6$ m to show local asymptotic stability for the closed-loop system (22) around a specific set of paths in backward motion. Then we show simulation results where we have perturbed the initial position of the system to analyze how well the LQ-controller handles disturbances by verifying a decreasing Lyapunov function. With the weight matrix $Q = \text{diag}([0.05, 10, 8, 2])$ and $v_3 < 0$ the feedback gain in the LQ-controller becomes

$$K = \begin{bmatrix} 0.22 & -4.88 & 6.18 & -3.84 \end{bmatrix} \quad (33)$$

In the experiments an eight shaped path is generated, e.g. by a motion planner as the one in [12], that satisfies the true system dynamics (1).

### A. Stability around a set of paths

To be able to apply Theorem 2 on the closed-loop system (22) in backward motion ($v_3 = -1$) we need to specify the set where the desired path lies in. From the Appendix we can conclude that the Jacobian matrix (27) only explicitly depends on the parameters $\beta_{3,0}$, $\beta_{2,0}$ and $u_0$. The parameters $z_{3,0}$ and $\tilde{\theta}_{3,0}$ do not vary because they have to be identically zero. The path parameters can be bounded in $\mathbb{P}$ and $\mathbb{U}$ by taking the maximum and the minimum values for each parameter independently along the eight shaped path. However, if doing so we neglect the fact that there is a natural connection between them. The simulated path (blue lines) can be seen in Fig. 3 where we also have outer bounded the path parameters (red lines) and thus define $\mathbb{M} = \mathbb{P} \times \mathbb{U} = \{(\beta_{3,0}, \beta_{2,0}, u_0) \in \mathbb{R}^3\}$ such that

$$|\beta_{3,0}| \leq 40 \text{ deg} \qquad |\beta_{2,0}| \leq 20 \text{ deg} \qquad (34)$$
$$|u_0| \leq 0.37 \qquad |\beta_{2,0} - \beta_{3,0}| \leq 20 \text{ deg}$$
$$|\arctan u_0 - \beta_{2,0}| \leq 10 \text{ deg}$$

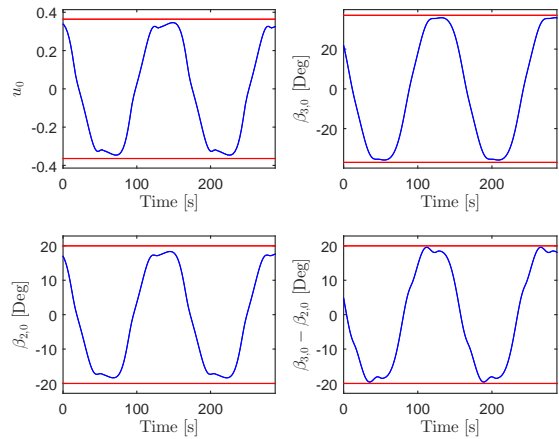

Fig. 3: The eight shaped reference path in blue together with four of the five constraints (34) in red. The time is not relevant for the analysis.

When the set $\mathbb{M}$ is specified we can solve the optimization problems (29a) and (29b) using the global solver bmibnb in YALMIP [15]. The obtained ranges for the elements in $A_{cl}$ are illustrated in Fig. 4. To show local asymptotic stability of the closed-loop system (22) for $(\mathbf{p}_{e,0}, u_0) \in \mathbb{M}$ we solve the LMIs (32) in Matlab using YALMIP [15]. For numerical reasons we bound $P$ from below and for robustness properties we introduce a scalar variable $\mu \geq 1$ and define the optimization problem as

$$\text{minimize} \quad \mu \qquad (35)$$
$$\text{subject to} \quad \mu I_{4\times 4} \succeq P \succeq I_{4\times 4}, \text{ and } (32)$$

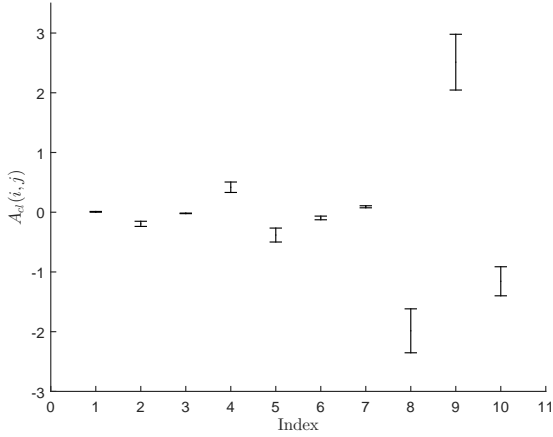

Fig. 4: A two dimensional illustration of the ten dimensional box describing the conservative approximation of the convex hull $\bar{A}_{cl}$ when $(\mathbf{p}_{e,0}, u_0) \in \mathbb{M}$. Each vertical line denotes the interval in one of the ten dimensions.

Putting the decay rate $\varepsilon = 0.001$ in (32) the optimal solution to (35) is $\mu = 118.14$ and the common symmetric positive definite matrix $P$ is:

$$P = \begin{bmatrix} 1.34 & -5.67 & 2.37 & -0.015 \\ -5.67 & 102.99 & -32.88 & 2.91 \\ 2.37 & -32.88 & 44.41 & -1.45 \\ -0.015 & 2.91 & -1.45 & 2.41 \end{bmatrix} \qquad (36)$$

### B. Simulation results

To evaluate how the LQ-controller handles disturbances we simulate the system with a perturbation in all initial states $\tilde{\mathbf{p}}_e(0) = [-4.2, -0.1, 0.1, -0.3]^T$. The velocity of the trailer is set to $-1$ m/s. The simulation results is provided in Fig. 5-8. In Fig. 5 the resulting path taken by the rear axle of the trailer is plotted together with the desired eight shaped path $(x_{3,0}, y_{3,0})(s)$. A similar plot can be viewed in Fig. 6 where all four error states $\tilde{\mathbf{p}}_e$ are plotted. One can conclude that the controller smoothly manages to stabilize the truck and trailer system around the nominal path even with quite big initial errors in all states. In Fig. 7 the feedback control signal $\tilde{u}$ is plotted over time. From this figure we see that the feedback part of the LQ controller mostly handles disturbance rejection and the feed forward part of the controller handles path tracking. We have already concluded that the closed-loop system (22) is local asymptotically stable if the system starts sufficiently close to the desired path. Using the Lyapunov matrix $P$ in (36) the quadratic Lyapunov function $V(t) = \tilde{\mathbf{p}}_e^T(t) P \tilde{\mathbf{p}}_e(t)$ together with its derivative $\dot{V}(t) = 2\tilde{\mathbf{p}}_e^T(t) P \dot{\tilde{\mathbf{p}}}_e$ during the simulation are shown in Fig. 8. This figure confirms the result by Theorem 2, i.e., that the nonlinear closed-loop system (22) is local asymptotically stable around paths with the path parameters bounded in $\mathbb{M}$ since $V(t) \to 0$ monotonically as $t \to \infty$.

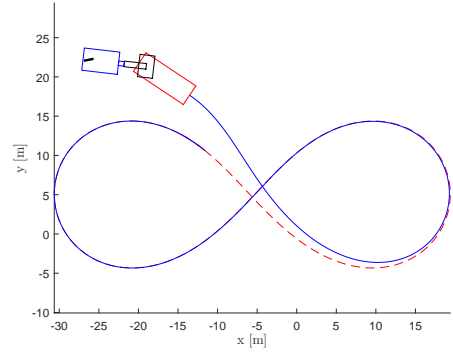

Fig. 5: Simulation result when tracking the eight shaped reference path. Blue line is the simulated states $(x_3, y_3)(t)$ and the red dotted line represents the reference path $(x_{3,0}, y_{3,0})(s)$.

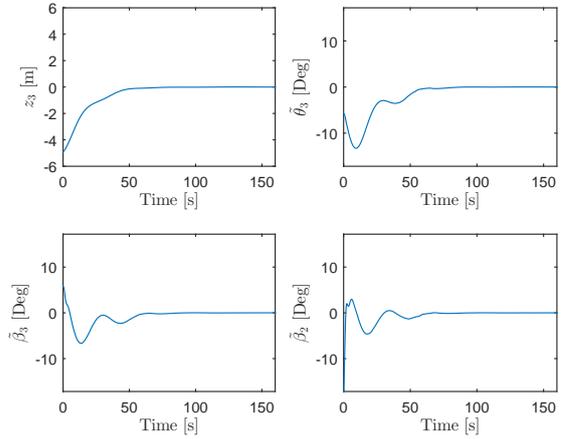

Fig. 6: Simulation result when tracking the eight shaped reference path. The blue lines are the error states where all states converge to zero.

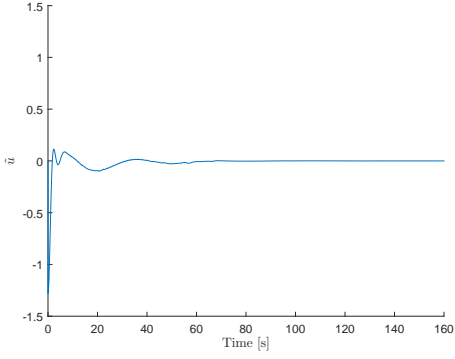

Fig. 7: Simulation result when tracking the eight shaped reference path. The blue line is the feedback control signal $\tilde{u}$.

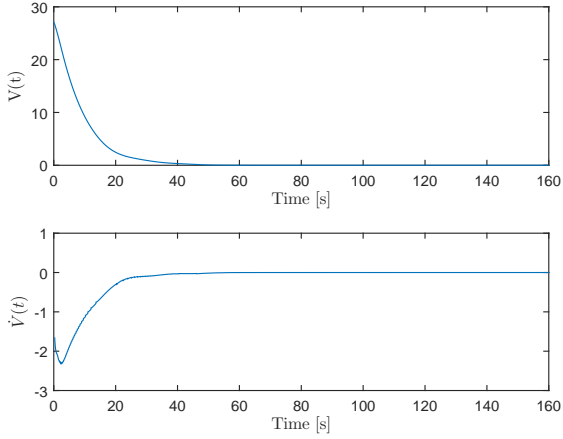

Fig. 8: The Lyapunov function $V(t) = \tilde{\mathbf{p}}_e^T P \tilde{\mathbf{p}}_e$ and its derivative $\dot{V}(t) = 2\tilde{\mathbf{p}}_e^T(t) P \dot{\tilde{\mathbf{p}}}_e$ during the simulation. As predicted, the closed-loop system is asymptotically stable since the Lyapunov function is monotonically decreasing to zero.

## VI. CONCLUSIONS AND FUTURE WORK

In this work a path following controller with local asymptotic stability guarantees is presented for a general 2-trailer system with possible off-axle hitching around a set of paths in backward motion. The result is obtained by combining techniques from global optimization with LMI techniques for LDIs with the aim of finding a Lyapunov function common to all paths in some predefined set, typically chosen to represent the possible output from a path planner. The theoretical result is confirmed in practice by simulations. As future work we would like to decrease the conservativeness of the result in order to be able to guarantee stability for a larger set of paths as well as for other relevant model uncertainties. Ongoing work includes implementing the controller and evaluating the performance on a full-scale truck and trailer system.

## VII. APPENDIX

The linearization of the nonlinear model (7b)-(7e) around origin yields,

$$\dot{\tilde{\mathbf{p}}}_e = A(\bar{0}, 0, \mathbf{p}_{e,0}, u_0)\tilde{\mathbf{p}}_e + B(\bar{0}, 0, \mathbf{p}_{e,0}, u_0)\tilde{u} \tag{37}$$

where

$$A(\tilde{\mathbf{p}}_e, \tilde{u}, \mathbf{p}_{e,0}, u_0) = \frac{\partial f_e}{\partial \tilde{\mathbf{p}}_e}(\tilde{\mathbf{p}}_e, \tilde{u}, \mathbf{p}_{e,0}, u_0) \tag{38}$$

$$B(\tilde{\mathbf{p}}_e, \tilde{u}, \mathbf{p}_{e,0}, u_0) = \frac{\partial f_e}{\partial \tilde{u}}(\tilde{\mathbf{p}}_e, \tilde{u}, \mathbf{p}_{e,0}, u_0) \tag{39}$$

The matrices $A(\tilde{\mathbf{p}}_e, \tilde{u}, \mathbf{p}_{e,0}, u_0)$ and $B(\tilde{\mathbf{p}}_e, \tilde{u}, \mathbf{p}_{e,0}, u_0)$ are some very complex functions and the structure can be written as

$$A = v_3 \begin{bmatrix} 0 & a_{12} & 0 & 0 \\ a_{21} & a_{22} & a_{23} & 0 \\ a_{31} & a_{32} & a_{33} & a_{34} \\ a_{41} & a_{42} & a_{43} & a_{44} \end{bmatrix} \quad B = v_3 \begin{bmatrix} 0 \\ 0 \\ b_3 \\ b_4 \end{bmatrix} \tag{40}$$

For compactness we use $\mathbf{p}_e = \tilde{\mathbf{p}}_e + \mathbf{p}_{e,0}$ and $u = \tilde{u} + u_0$, then the expressions for each position of these matrices can be written as

$$a_{12} = \cos \tilde{\theta}_3 \tag{41}$$

$$a_{21} = -\frac{\cos \tilde{\theta}_3 \tan^2 \beta_{3,0}}{(L_3 - z_3 \tan \beta_{3,0})^2} \tag{42}$$

$$a_{22} = \frac{\sin \tilde{\theta}_3 \tan \beta_{3,0}}{L_3 - z_3 \tan \beta_{3,0}} \tag{43}$$

$$a_{23} = \frac{1 + \tan^2 \beta_3}{L_3} \tag{44}$$

$$a_{31} = -\frac{L_3 \cos \tilde{\theta}_3 \tan \beta_{3,0}}{(L_3 - z_3 \tan \beta_{3,0})^2} \left( \frac{\tan \beta_{2,0} - \frac{M_1}{L_1} u_0}{L_2 \cos \beta_{3,0} C_2(\beta_{2,0}, u_0)} - \frac{\tan \beta_{3,0}}{L_3} \right) \tag{45}$$

$$a_{32} = \frac{L_3 \sin \tilde{\theta}_3}{L_3 - z_3 \tan \beta_{3,0}} \left( \frac{\tan \beta_{2,0} - \frac{M_1}{L_1} u_0}{L_2 \cos \beta_{3,0} C_2(\beta_{2,0}, u_0)} - \frac{\tan \beta_{3,0}}{L_3} \right) \tag{46}$$

$$a_{33} = \frac{\sin \beta_3 \left( \tan \beta_2 - \frac{M_1}{L_1} u \right)}{L_2 \cos^2 \beta_3 C_2(\beta_2, u)} - \frac{1 + \tan^2 \beta_3}{L_3} \tag{47}$$

$$a_{34} = \frac{(1 + \tan^2 \beta_2)}{L_2 C_2(\beta_2, u) \cos \beta_3} \left( 1 + \frac{\left( \frac{M_1}{L_1} u - \tan \beta_2 \right) \frac{M_1}{L_1} u}{C_2(\beta_2, u)} \right) \tag{48}$$

$$a_{41} = -\frac{L_3 \cos \tilde{\theta}_3 \tan \beta_{3,0}}{(L_3 - z_3 \tan \beta_{3,0})^2} \left( \frac{\frac{u_0}{L_1 \cos \beta_{2,0}} - \frac{\tan \beta_{2,0}}{L_2} + \frac{M_1}{L_1 L_2} u_0}{\cos \beta_{3,0} C_2(\beta_{2,0}, u_0)} \right) \tag{49}$$

$$a_{42} = \frac{L_3 \sin \tilde{\theta}_3}{L_3 - z_3 \tan \beta_{3,0}} \left( \frac{\frac{u_0}{L_1 \cos \beta_{2,0}} - \frac{\tan \beta_{2,0}}{L_2} + \frac{M_1}{L_1 L_2} u_0}{\cos \beta_{3,0} C_2(\beta_{2,0}, u_0)} \right) \tag{50}$$

$$a_{43} = \frac{\sin \beta_3 \left( \frac{u}{L_1 \cos \beta_2} - \frac{\tan \beta_2}{L_2} + \frac{M_1}{L_1 L_2} u \right)}{\cos^2 \beta_3 C_2(\beta_2, u)} \tag{51}$$

$$a_{44} = \frac{\left( \frac{u \sin \beta_2}{L_1 \cos^2 \beta_2} - \frac{1 + \tan^2 \beta_2}{L_2} \right)}{\cos \beta_3 C_2(\beta_2, u)} \\ - \frac{M_1}{L_1} u(1 + \tan^2 \beta_2) \frac{\left( \frac{u}{L_1 \cos \beta_2} - \frac{\tan \beta_2}{L_2} + \frac{M_1}{L_1 L_2} u \right)}{\cos \beta_3 C_2^2(\beta_2, u)} \tag{52}$$

and

$$b_3 = -\left(\frac{M_1}{L_1 L_2 \cos\beta_3 C_2(\beta_2, u)} + \frac{\left(\tan\beta_2 - \frac{M_1}{L_1}u\right)\frac{M_1}{L_1}\tan\beta_2}{L_2 \cos\beta_3 C_2^2(\beta_2, u)}\right) \tag{53}$$

$$b_4 = \frac{1}{\cos\beta_3 C_2(\beta_2, u)} \left(\frac{1}{L_1 \cos\beta_2} + \frac{M_1}{L_1 L_2}\right.$$
$$\left. - \frac{M_1 \tan\beta_2}{L_1 C_2(\beta_2, u)} \left(\frac{u}{L_1 \cos\beta_2} - \frac{\tan\beta_2}{L_2} + \frac{M_1}{L_1 L_2}u\right)\right) \tag{54}$$

## REFERENCES

[1] M. Werling, P. Reinisch, M. Heidingsfeld, and K. Gresser, "Reversing the general one-trailer system: Asymptotic curvature stabilization and path tracking," *Intelligent Transportation Systems, IEEE Transactions on*, vol. 15, no. 2, pp. 627–636, 2014.
[2] O. J. Sordalen, "Conversion of the kinematics of a car with n trailers into a chained form," *IEEE J. Robot. Automat.*, 1993.
[3] P. Rouchon, M. Fliess, J. Lévine, and P. Martin, "Flatness, motion planning and trailer systems," in *Decision and Control, 1993., Proceedings of the 32nd IEEE Conference on*. IEEE, 1993, pp. 2700–2705.
[4] M. Sampei, T. Tamura, T. Kobayashi, and N. Shibui., "Arbitrary path tracking control of articulated vehicles using nonlinear control theory," *Control Systems Technology, IEEE Transactions on*, 1995.
[5] D. Tilbury, R. M. Murray, and S. Shankar Sastry, "Trajectory generation for the n-trailer problem using goursat normal form," *Automatic Control, IEEE Transactions on*, vol. 40, no. 5, pp. 802–819, 1995.
[6] C. Altafini, "The general n-trailer problem: conversion into chained form," in *Proc. IEEE Decision and Control*, 1998.
[7] P. Bolzern, R. M. DeSantis, A. Locatelli, and D. Masciocchi, "Path-tracking for articulated vehicles with off-axle hitching," *Control Systems Technology, IEEE Transactions on*, vol. 6, no. 4, pp. 515–523, 1998.
[8] C. Altafini, "Path following with reduced off-tracking for multibody wheeled vehicles," *Control Systems Technology, IEEE Transactions on*, vol. 11, no. 4, pp. 598–605, 2003.
[9] C. Altafini, A. Speranzon, and K. H. Johansson, *Hybrid Control of a Truck and Trailer Vehicle*. Springer Berlin Heidelberg, 2002.
[10] N. Evestedt, O. Ljungqvist, and D. Axehill, "Path tracking and stabilization for a reversing general 2-trailer configuration using a cascaded control approach." arXiv:1602.06675, 2016.
[11] J. David and P. Manivannan, "Control of truck-trailer mobile robots: a survey," *Intelligent Service Robotics*, vol. 7, no. 4, pp. 245–258, 2014.
[12] O. Ljungqvist, "Motion planning and stabilization for a reversing truck and trailer system," 2016.
[13] S. P. Boyd, L. El Ghaoui, E. Feron, and V. Balakrishnan, *Linear matrix inequalities in system and control theory*. SIAM, 1994, vol. 15.
[14] M. Pakmehr, N. Fitzgerald, E. M. Feron, J. S. Shamma, and A. Behbahani, "Gain scheduling control of gas turbine engines: stability by computing a single quadratic lyapunov function," in *ASME Turbo Expo 2013: Turbine Technical Conference and Exposition*. American Society of Mechanical Engineers, 2013, pp. V004T06A027–V004T06A027.
[15] J. Löfberg, "Yalmip: A toolbox for modeling and optimization in matlab," in *Computer Aided Control Systems Design, 2004 IEEE International Symposium on*. IEEE, 2004, pp. 284–289.
[16] M. Sampei and K. Furuta, "On time scaling for nonlinear systems: Application to linearization," *Automatic Control, IEEE Transactions on*, vol. 31, no. 5, pp. 459–462, 1986.
[17] B. D. Anderson and J. B. Moore, *Optimal control: linear quadratic methods*. Courier Corporation, 2007.
[18] H. K. Khalil and J. Grizzle, *Nonlinear systems*. Prentice hall New Jersey, 1996, vol. 3.